\documentclass[12pt]{amsart}
\usepackage{bookman}
\usepackage[pagebackref]{hyperref}
\usepackage[alphabetic,backrefs,lite,nobysame]{amsrefs} % for bibliography
\usepackage{amsmath}
\usepackage{url,amssymb,enumerate,colonequals,bm}
\usepackage[headings]{fullpage}
\usepackage[all]{xy}
\usepackage{comment}
\numberwithin{equation}{section}
\numberwithin{figure}{section}

\usepackage{mathrsfs}% for \mathscr (script letters)
\usepackage{tikz-cd}

\usepackage{calrsfs}
\newcommand{\ab}{{\operatorname{ab}}}

\newcommand{\HH}{{\operatorname{H}}}

\newcommand{\calL}{{\mathcal L}}
\newcommand{\calA}{{\mathcal A}}

\newcommand{\Z}{{\mathbb Z}}

\newcommand{\C}{{\mathbb C}}

\newcommand{\calB}{{\mathcal B}}

\newcommand{\bL}{\mathbf L}

\newtheorem*{theorem*}{Theorem}
\newtheorem*{conjecture*}{Conjecture}
\newtheorem{conjecture}[equation]{Conjecture}
\newtheorem*{question*}{Question}

\newtheorem{theorem}[equation]{Theorem}

\newtheorem{proposition}[equation]{Proposition}

\theoremstyle{definition}
\newtheorem{definition}[equation]{Definition}

\newtheorem{remark}[equation]{Remark}

\theoremstyle{remark}

\newcommand{\oo}{\mathcal{O}}
\newcommand{\Q}{\mathbb{Q}}
\newcommand{\R}{\mathbb{R}}

\newcommand{\Gal}{\mathrm{Gal}}

\def\vv{\vskip10pt}

\newtheorem*{lemma*}{Lemma}
\newtheorem*{notation*}{Notation}
\newtheorem*{proposition*}{Proposition}

\title[The ring of integers of any abelian extension of $\Q$ is undecidable]{The first-order theory of the ring of integers of any abelian extension of $\Q$ is undecidable}
\author{Barry Mazur}
\address{Department of Mathematics\\ Harvard University \\ Cambridge, MA 02138-2901}
\email{mazur@g.harvard.edu}
\urladdr{http://www.math.harvard.edu/~mazur}
\author{Karl Rubin}
\address{Department of Mathematics\\
UC Irvine\\
Irvine, CA 92697, 
USA}
\email{krubin@uci.edu}
\urladdr{https://math.uci.edu/~krubin}
\author{Alexandra Shlapentokh}
\address{Department of Mathematics \\ East Carolina University \\ Greenville, NC 27858}
\email{shlapentokha@ecu.edu }
\urladdr{myweb.ecu.edu/shlapentokha}
\begin{document}

\begin{abstract}
  We  show that any ring of totally real abelian integers has a JR-number and therefore the first-order theory of every ring of totally real abelian integers is undecidable. 
It then follows  that the first-order theory of the ring of integers of any abelian extension of $\Q$ is undecidable and that there is one formula without parameters defining $\Z$ in any such ring. The main arithmetic input is a recent paper of Bajpai, Das, Kedlaya, and Mello proving that the houses of cyclotomic integers
form a well-ordered subset of $\R$.

\end{abstract}
\thanks{We acknowledge AI help in finding the relevant papers on the subject and proofreading our paper. The research for this paper was partially supported by DMS grants 2152098 (AS), 2152149 (BM), 2152262 (KR)}

\maketitle
%\makeatother
\tableofcontents
\addtocontents{toc}{\setcounter{tocdepth}{2}} 

\section{Our main results}
\vv
We answer a long-standing question concerning the elementary theory of the ring of integers of {\it any} abelian extension of $\Q$. 
\begin{theorem}
\label{thm:main}
The first-order theory of the ring of integers of {\it any}\footnote{A recent result of C. Springer \cite{Sp26} shows undecidability of the first-order theory of the rings of integers of {\it a certain class of abelian extensions of $\Q$} ({\it including the ring  $\Z^{ab}$ of integers in  $\Q^{ab}$; i.e., the largest abelian extension of $\Q$}).} abelian extension of $\Q$ is undecidable. 
\end{theorem}

We have previously shown in \cites{MRSh24-1}---see Section \ref{1} below--- that to prove Theorem \ref{thm:main} we need only consider totally real abelian extensions. For these extensions we show---see Section \ref{pf} below---  the following.
\begin{proposition}
\label{thm:classify}
 Any ring of totally real abelian integers has a JR-number (Julia Robinson number, see Definition \ref{def:JR}).
 \end{proposition}
J. Robinson showed that {\it any} totally real ring $R$ of algebraic integers with a JR-number is finitary (see Definition \ref{def:fin}) and there exists a first-order formula $\psi(x,\bar t)$, where $x$ is the only variable not in the range of any quantifier, such that $R\models \psi(x,\bar t)$ \footnote{$R\models \phi(x,\bar t)$ means $\phi(x,\bar t)$ is true over $R$.} if and only if $x \in \Z$ (\cite{Rob3}).   In other words, there exists  a uniform  definition of $\Z$ without parameters across all such rings $R$.  Thus, we get:
\begin{proposition}\label{prop} \begin{enumerate} \item There is a first order formula without parameters that  defines  $\Z$ over {\rm\bf any subring} of the ring of totally real abelian integers. \vv
\item 
There is a first order formula without parameters that  defines  $\Z$ over {\rm\bf the ring of integers} of any  abelian extension of $\Q$. \end{enumerate}
\end{proposition}
Theorem \ref{thm:main} follows from  Proposition \ref{prop} (2).\\

\section{Background Terminology}
\label{sec:lang}
 The first-order language of rings is the language $\calL_R=(0, 1, +, \times)$.
 \begin{itemize}\item By a {\it first-order formula} $\phi(y_1,...,y_s)$ is equivalent to
\[
    Q_1 x_1 ... Q_r x_r \theta(y_1,...,y_s,x_1,...,x_r),
\]

  where each $Q_i$ is existential or universal and $\theta$ is a Boolean
  combination of polynomial equations with integer coefficients.\vv
 \item Call it  a {\it Diophantine formula}  if all the $Q_i$ are  existential.
 \end{itemize} \vv

Given a ring $R$, we can define a subset $A$ of $R$ in the following manner: 
\[
A=\{y|R\models \phi(y,\bar x)\},
\]
where $\phi(y,\bar x)$ is a first order formula with $y$ not being in the range of any quantifier and all variables constituting $\bar x$ being in the range of some quantifier in the formula.\\
 
The following notions play a key role in the results of this paper.
\begin{definition}[Finitary collections]
\label{def:fin}
We will call a collection of sets {\bf finitary} if it contains finite sets of arbitrarily large size.  
\end{definition}
\begin{definition}[The house of an algebraic integer]
For an algebraic integer $\alpha$ define its house by
\[
\HH(\alpha)=\max_{\sigma:\Q(\alpha)\hookrightarrow\mathbb C}
                  |\sigma(\alpha)|.
\]
\end{definition}
\begin{definition}[JR-numbers]
\label{def:JR}
Let $R$ be a ring of totally real integers.  A JR-number for the ring $R$ (denoted by JR($R$)) is either a real number or $\infty$ satisfying the following conditions.
\begin{itemize}
\item A positive real number $t$ is JR($R$) if the interval $[0,t]$ contains infinitely many houses of totally positive elements of $R$ and for any $0<\varepsilon<t$ the interval $[0,t-\varepsilon]$ contains finitely many houses of totally positive elements of $R$.
\item JR($R$)$=\infty$ if for every positive real number $s$ the interval $[0,s]$ contains finitely many houses of totally positive elements of $R$.
\end{itemize}
\end{definition}
\begin{remark}
If $x, y$ are totally real algebraic integers and $\HH(x)=\HH(y)$, then either $x$ and $y$ are conjugate over $\Q$ or $-x$ and $y$ are conjugate over $\Q$.  Therefore, for every real number $h$ the set of algebraic integers $x$ with  $\HH(x)=h$ is finite or empty.
\end{remark}
\begin{definition}[Finitary rings]
Let $R$ be a ring.  Suppose there exists a first-order formula $\psi(\bar a, x)$ in the language of rings such that the collection of sets $\calB$ consisting of the sets $B_{\bar a}=\{x \in R|R\models \psi(\bar a,x)\}$, where $\bar a$ ranges over all tuples of appropriate arity with elements in $R$, is finitary.  Then we call $R$ a finitary ring.
\end{definition}\vv

\section{ Passage from  totally real  abelian algebraic extensions of $\Q$ to arbitrary  abelian algebraic extensions of $\Q$}\label{1}

   In \cite{MRSh24-1}[Proposition 2.11] we proved the following proposition.
   \begin{proposition}
   Let $K \ne \Q$ be a totally real extension of $\Q$ and let $L$ be an extension of degree 2 of $K$. Let $\hat R \subseteq \oo_L$ be any subring of $\oo_L$ not contained in the ring of integers of a quadratic imaginary extension of $\Q$. Then there exists a ring of totally real integers $R \subset \hat R$ and a first-order formula $\psi(z)$ without parameters where $z$ is the only free variable such that for all $z\in \hat R$ we have that $\hat R\models \psi(z)$ if and only if $z \in R$. Further, this formula is not dependent on $L$ or $\hat R$.  
   \end{proposition}
 \begin{proposition}
\label{prop:real}
If there exists a uniform definition of $\Z$ without parameters over every ring of totally real abelian algebraic integers, then there  exists a uniform definition of $\Z$ without parameters over the ring of integers of every abelian extension of $\Q$.
\end{proposition}
\begin{proof}
First, we note that every non-totally-real abelian extension of $\Q$ is an extension of degree 2 of a totally real field. For the sake of convenience, we embed $\bar \Q$, the algebraic closure of $\Q$, in $\C$.   Now, let ${\bf K}$ be an abelian, possibly infinite, extension of $\Q$.  Complex conjugation $\sigma$ is an element of $\Gal({\bf K}/\Q)$.  Let ${\bf K}^{\R} \subset \R$ be the fixed field of $\sigma$.  We claim ${\bf K}^{\R}$ is totally real.  In other words, for any $\tau \in \Gal({\bf K}/\Q)$ we have that $\tau({\bf K}^{\R}) \subset \R$.  Indeed, for any $x \in {\bf K}^{\R}$ we have that $\sigma(\tau(x))=\tau(\sigma(x))=\tau(x)$.  Since $\tau(x)$ is in the fixed field of complex conjugation we have that $\tau(x) \in \R$.  

Second, we note that in \cite{MRSh24-1} we constructed a uniform existential definition of $\Z$ without parameters over every imaginary quadratic extension of $\Q$.  Third, to distinguish an imaginary quadratic extension from other extensions of degree 2 of totally real fields it is enough to see if the fields have at least seven units.  Indeed, there are only six units that can show up in an imaginary quadratic extension: $\pm 1, \pm i, \pm \zeta_3 \text{ (a third root of unity)}$.  So, if $\phi(x)$ is a first-order formula without parameters defining $\Z$ over the ring of integers of any non-trivial totally real extension of $\Q$ and any  extension of degree 2 of a totally real field $K\ne \Q$ and $\theta(x)$ is a first-order formula without parameters defining $\Z$ over the ring of integers of any imaginary quadratic extension, then the following formula without parameters will define $\Z$ over the rings of integers of any abelian extension $L$ of $\Q$. 
\[
\Gamma(x)=
\]
\[
[\exists \mu_1,\ldots, \mu_7, \nu_1,\ldots, \nu_7 \in \oo_L: \mu_i \ne \mu_j \text{ for } i\ne j \land \mu_i\nu_i=1 \land \phi(x)] \lor 
\]
\[
[\forall \mu_1,\ldots,\mu_7 \in \oo_L [(\exists i\ne j: \mu_i = \mu_j) \lor (\exists i \forall u \in \oo_L: u \mu_i \ne 1)] \land \theta(x)].
\]

\section{Proof of our main theorem (Theorem \ref{thm:classify})}\label{pf}
 We will be using the following result\footnote{ We thank AI for bringing this result to our attention.}:\vv
\begin{proposition}[Bajpai--Das--Kedlaya--Mello]\label{thm:gap}
For each $C>0$ there is $\varepsilon>0$ such that no cyclotomic integer
$\alpha$ satisfies
\begin{equation}
\label{eq:without}
 C<\HH(\alpha)<C+\varepsilon.
\end{equation}\end{proposition}
\vv
Note that 
for $R$  a ring of totally real abelian integers  the set of houses $\HH(R)$ of all totally positive  elements of $R$ is a well-ordered set. We consider two possibilities:\vv
\begin{itemize} \item  $\HH(R)$  does not have any accumulation points, in which case JR($R$) = $\infty$,
\vv
or\vv \item $\HH(R)$ has an accumulation point. 
\end{itemize}\vv 
 
 Suppose $\HH(R)$ has accumulation points.  Denote the set of these points by $A(R)$.  We claim $A(R)$ has a minimum.  Suppose not.  Note that $A(R)$ is bounded from below.  Therefore, it has a greatest lower bound $L$.  If it is not an element of the set, then for any positive $\varepsilon$ there are infinitely many elements $a$ of $A(R)$ such that $a-L<{\frac{\varepsilon}{2}}$.  Since $a$ is an accumulation point of $H(R)$ there are infinitely many totally positive elements $x$ of $R$ such that $|a-H(x)|<{\frac{\varepsilon}{2}}$.  Observe that by construction $H(x)>L$ and 
 \[
 H(x)-L=H(x)-a+a-L<|H(x)-a|+|a-L|< \varepsilon.
 \]
 This contradicts \eqref{eq:without}.  Hence $A(R)$ has a minimum.  
 
 Let $A$ be the smallest accumulation point.  Then by \eqref{eq:without} for some $\varepsilon$ we have that the interval $(A, A+\varepsilon)$ contains no houses.  Let $\delta <\varepsilon$ and consider the set of all houses $h$ such that $|h-A|<\delta$.  Then all such houses must be in the interval $[0,A]$ and there must be infinitely many of them.  Now let $B \in (0,A)$.  Then $B$ is not an accumulation point of $\HH(R)$ and, by the Bolzano-Weierstrass theorem, there can be only finitely many houses in the interval $[0,B)$.  
 If $\HH(R)$ has no accumulation points, then for any $B>0$, again by the Bolzano-Weierstrass theorem,  the  interval $(0,B)$ has finitely many points, so JR($R$) = $A$.
 
Finally, if the ring $R$ has no accumulation points, then every finite interval contains at most finitely many houses.  Therefore, in this case $\text{JR}(R)=\infty$. 
 \end{proof}

There have been a number of papers investigating JR-numbers and the elementary theory of rings of integers of totally real extensions of $\Q$.  See, for example, \cites{V3}, \cites{VidVid1}, \cites{VV16}, \cites{GilRan}, \cites{CasVidVid}, \cites{MarUtrVid}, \cites{Sp}, \cites{CF21},  \cites{Sp23}.

\section{Appendix: Existential definability results}
In some cases one can show that the positive {\bf existential} theory of infinite abelian extensions is undecidable.
\subsection{Rings of integers of infinite extensions of $\Q$ with undecidable existential theory}
There exist abelian extensions of $\Q$ where we know the existential theory of the rings of integers to be undecidable.  The following result is due to  
the third author \cites{Sh37},  K. Kato \cites{Kato}, K. Ribet \cites{Rib} 
and D. Rohrlich \cites{Rohr1, Rohr2} (see \cites{LR}[Theorem 1.2]).
\begin{theorem}\label{KRR}
Let $\bL$ be an abelian extension with finitely many ramified primes.  
Then $\Z$ is existentially definable over $\oo_{\bL}$.
\end{theorem}
This result was one of the early examples of using diophantine stability (in this case of an elliptic curve) to show existential definability of $\Z$ over infinite algebraic extensions of $\Q$.  A more general approach was developed in \cite{MRSh24} where the following theorem was proved. 
\begin{theorem}
\label{pre-mainth1}
Let $K$ be a number field, let $L$ be an algebraic  extension of $K$, possibly infinite degree, and let
$\oo_K \subset \oo_L$ be their respective rings of integers.  Suppose $A$ is an abelian variety 
defined over $K$ such that $A(K)$ is infinite and $A(L)/A(K)$ is a torsion group.  
If at least one of the following conditions is satisfied:
\begin{enumerate}
\item
$L$ is a number field,
\item
$L$ is totally real,
\item
$L$ is a quadratic extension of a totally real field,
\end{enumerate}
then $\oo_K$ has a diophantine definition over $\oo_L$. 
\end{theorem}

Using some heuristics for elliptic curves the authors conjectured the following in \cite{MRSh24}.
\begin{conjecture}
\label{CON1}
For the primes $p=7$, $11$, or $13$ there are subfields 
${\bL} \subset {\Q}^{\ab}$ for which the field extension 
$ {\Q}^{\ab}/{\bL}$ is cyclic of degree $p$ and such that 
$\oo_{\bL}$, the ring of integers in ${\bL}$, is diophantine 
undecidable. 
\end{conjecture}
In fact, we conjecture (the stronger statement) that ${\Z}$ is existentially 
definable over ${\mathcal O}_{\bL}$ for the fields ${\bL}$ in Conjecture \ref{CON1}. 
We do not make the same conjecture for $\oo_{{\Q}^{\rm ab}}$.

There have been many results and conjectures concerning 
existential definability of $\Z$ over the ring of integers of number fields and diophantine stability of elliptic curves. In particular, we now know that the existential theory of the ring of integers of any number field is undecidable (see \cite{KP24}, \cite{ABHS25}, \cite{Z25}).

\section{Appendix: A uniform first-order definition of $\Z$ for ``non-big fields''}
      In \cite{MRSh24}  we constructed a definition $\Z$ over the rings of integers of some infinite algebraic extensions of $\Q$.  Whether this definition works for a particular field is not dependent on whether the ring in question is real or not.  Rather it depends on whether the field in question is ``big'' or not (see \cite{MR20}).

\begin{definition}
\label{big} 
Let $K$ be a field of algebraic numbers.  
We say that $K$ is {\em  big } if 
$$
[K :{\mathbb  Q}] = \prod_{p\ {\rm prime}} p^\infty.
$$  
Equivalently,
$K$ is big if for every positive integer $n$, $K$ contains a number field $F$ 
with $[F : {\Q}]$ divisible by $n$.
\end{definition}

In \cite{MRSh24-1} we proved the following.

\begin{theorem}
\label{thm:units}
Let $\calA$ be the collection of all non-big fields of algebraic numbers.
There exists a first-order formula of the form 
``$\forall \forall\exists \ldots \exists$'' uniformly defining $\Z$ over the 
rings of integers of all fields in $\calA$.  
In particular the first-order theory of these rings can be uniformly shown to be first-order undecidable.
\end{theorem}

\end{document}